\documentclass[12pt,a4paper]{amsart}

\input prepictex.tex
\input pictex.tex
\input postpictex.tex
\chardef\bslash=`\\ 

\makeatletter
\def\verbatim{\interlinepenalty\@M \@verbatim
  \leftskip\@totalleftmargin\advance\leftskip2pc
  \frenchspacing\@vobeyspaces \@xverbatim}
\makeatother
\newtheorem{theorem}{Theorem}[section]

\newtheorem{proposition}[theorem]{Proposition}

\newtheorem{remark}[theorem]{Remark}

\theoremstyle{definition}
\newtheorem{definition}{Definition}[theorem]

\numberwithin{equation}{section}

\newcounter{picture}

\newcommand{\RR}{{\Bbb R}}
\newcommand{\ZZ}{{\Bbb Z}}

\newcommand{\G}{{\Gamma}}
\newcommand{\Om}{{\Omega}}

\newsymbol\rtimes 226F 

\begin{document}

\title[]{ Negative definite kernels and a
 dynamical characterization of property (T)
 for countable groups }

\date{November 7, 1996}
\author[G. Robertson]{Guyan Robertson}
\author[T. Steger]{Tim Steger}
\address{Department of Mathematics  \\
        University of Newcastle\\  NSW  2308\\ AUSTRALIA}
\email{guyan@@maths.newcastle.edu.au }
\address{Istituto di Matematica e Fisica \\
Universit\`a di Sassari \\
via Vienna 2 \\ 07100 Sassari \\ ITALIA}
\email{steger@@ssmain.uniss.it }
\thanks{This research was supported by the Australian Research Council.} 
\thanks{ \hfill Typeset by  \AmS-\LaTeX}

\maketitle
\begin{abstract}
A  class of negative definite kernels is defined in terms of measure 
spaces. Using this concept, property (T) for a countable group $\G$ is
 characterized in terms of measure preserving actions
of $\G$, as follows.
If a set $S$ is translated a finite amount by
any fixed element of $\G$, then there is a
 uniform bound on how far $S$ is translated.
\end{abstract}

\section*{Introduction}

 A group $G$ has Kazhdan's property (T) if the trivial representation is
isolated in the
unitary dual of $G$ \cite{k}. We refer to \cite{hav} for an exposition of
many of the remarkable
properties of such groups.

Groups with property (T) are characterized by the absence of unbounded
negative definite functions \cite[Chapter~5, Theorem~20]{hav}. This
characterization 
was obtained in  \cite{del}, \cite{gui} and, independently, in \cite{akw}.
We use this
to prove the following dynamical characterization of property (T),
which is the main result of this paper. All measure spaces are positive
but not necessarily finite.
 
\medskip

{\sl
 A countable group $\G$ has property (T)
if and only if it satisfies the following condition:
for every measure-preserving action of $\G$ on a measure
space  $(\Om,{\mathcal B},\mu)$ and every set  $S \in {\mathcal B}$  such that
$\mu (S \triangle gS) < \infty $  for all  $g \in \G$, we have  
 $$\sup_{g\in \G} \mu(S \triangle gS) < \infty .$$
}

\begin{remark} \label{-1}
{\em
To see how the condition above can fail to hold,
suppose that there exists a homomorphism from $\G$ onto $\ZZ$.
(There is no such homomorphism, if $\G$ has property (T).)
Take  $\Omega = \ZZ$, $\mu =$ counting measure, $S = \{1,2,3,\dots\}$ and
let $\G$ act on $\Omega$
by translation. Then $\mu(S \triangle gS) < \infty $ for all $g \in \G$ but
$\sup_{g \in \G} \mu (S \triangle gS) = \infty $.
}
\end{remark}

\begin{definition}\label{pdk} Let $\G$ be a countable or finite set. Say
that a function 
$f : \G \times \G \longrightarrow [0, \infty)$ is a {\it measure definite
kernel} if there is a measure space $(\Omega,{\mathcal B}, \mu)$ which contains
 sets $S_x \in {\mathcal B}$ for $x \in \G$ such that $f(x,y) = \mu(S_x
\triangle S_y)$.
\end{definition}

\begin{remark} \label{0}
{\em
If $f$ is measure definite we can find 
a positive measure space $(\Omega,\mu)$ which contains
 sets $S_x$ for $x \in \G$ such that 
$f(x,y) = \mu(S_x \backslash S_y) = \mu(S_y \backslash S_x)$.
For this, simply replace $\Omega$ by $\Omega' = \Omega \times \{0,1\}$
with the product measure and replace each set $S_x$
by  $S_x' = (S_x \times \{0\}) \sqcup (S_x^c \times \{1\})$.
}
\end{remark}

\section{measure definite kernels}

We refer to \cite[Chapter~5]{hav} for the definition and properties 
of negative definite kernels. A negative definite kernel $f$
is assumed to be {\it normalized}, i.e. $f(x,x) = 0$.
A universal construction of a negative definite kernel
on a set $\G$ is to take a mapping
$x \mapsto v_x$ of $\G$ into a real Hilbert space  ${\mathcal H}$ and
to let the kernel $f(x,y) = \|v_x - v_y\|^2$.
Moreover, for a given $f$, the embedding $x \mapsto v_x$ of $\G$ in
${\mathcal H}$ is unique up to rigid motions of the closed convex hull of
$\{v_x : x \in \G \}$ \cite[Chapter~5, Proposition~14]{hav}.

\begin{proposition} \label{1}
Every measure definite kernel is negative definite.
\end{proposition}

{\sc Proof:} Let $f$ be a measure definite kernel. Using the notation of
Definition \ref{pdk},
fix $x_0 \in \G$ and define a function  $\eta :\G \to L^2(\Omega,\mu )$\ by
$\eta(x) = \chi_x - \chi_{x_0}$, 
where  $\chi_x$  denotes  the characteristic  function of ${S_x}$.
Then 
$$\| \eta(x) - \eta(y) \|^2 = \mu(S_x \triangle S_y) = f(x,y)$$
and so $f$ is negative definite \cite[Chapter~5.13, Exemples]{hav}.
\qed

\medskip
Every measure definite kernel is a pseudometric.
The negative definite kernel defined on $\RR$ by
$f(x,y) = (x - y)^2$ is clearly not a pseudometric.
The converse of Proposition~\ref{1} is therefore false.
However we shall show below that the
square root of any negative definite kernel is measure definite.
\medskip

For the rest of this section let $\G$ be a fixed countable set.
Let $\Omega_{00} = \{0,1\}^{\G}$ with the product topology.
Then $\Omega_{00}$ is either a Cantor set or finite.
Let $\Omega_0 = \Omega_{00} \backslash \{(0,0,0, \dots ), (1,1,1, \dots
)\}$
and let ${\mathcal B}_0$ be the $\sigma$-algebra of Borel sets in $\Omega_0$.
Given $x \in \G$, let $S_x^0 = \{ c \in \Omega_0 : c(x) = 1 \}$.
Our next result shows that ${\mathcal B}_0$ is a universal $\sigma$-algebra
for measure definite kernels on $\G$. 

\begin{proposition} \label{1+}
If $f : \G \times \G \longrightarrow [0, \infty)$ is a
measure definite kernel then
there exists a regular Borel measure $\mu_0$ on $(\Omega_0,{\mathcal B}_0)$
such that $f(x,y) =  \mu_0(S_x^0 \triangle S_y^0)$
for all $x,y \in \G$.
\end{proposition}

{\sc Proof:}
Suppose that a measure space $(\Omega, {\mathcal B}, \mu)$ is given, which
contains
 sets $S_x$,  such that  $\mu(S_x \triangle S_y)$
is finite for $x,y \in \G$. We can and will delete from $\Omega$ all points belonging to all
(respectively none) of the sets $S_x$ without changing $\mu(S_x \triangle
S_y)$.
Define a map $\phi : \Omega \to \Omega_0$ by
$$(\phi(p))(x)  =  \left\{ \begin{array}{ll}
                0 & \mbox{if} \quad p \notin S_x\\
                1  & \mbox{if} \quad p \in S_x .
                \end{array}
\right.$$

\noindent Since $\phi$ is measurable, we may define 
a corresponding measure 
$\mu_0 = \phi_*(\mu)$ on $(\Omega_0, {\mathcal B}_0)$.

Then 
$$
\mu_0(S_x^0 \triangle S_y^0)
= \mu(\phi^{-1}(S_x^0 \triangle S_y^0))
= \mu(S_x \triangle S_y) .
$$
\noindent In this way any measure definite kernel may be expressed as 
$\mu_0(S_x^0 \triangle S_y^0)$ for some measure $\mu_0$ on $(\Omega_0,{\mathcal
B}_0)$.
We claim (a)  every open set in $\Omega_0$ is $\sigma$-compact ;
(b)  $\mu_0(K) < \infty$ for every compact set $K$ in $\Omega_0$. 
By \cite[Theorem 2.18]{r}, these two assertions imply
that the measure $\mu_0$ is regular.

Let $\mathcal C$ denote the family of all cylinder sets of the
form  
\begin{equation*} 
\begin{split}
A
&=  (S_{x_1}^0\cap \dots \cap S_{x_m}^0) \backslash
(S_{y_1}^0\cup \dots \cup S_{y_n}^0)\\
&= \{ c \in \Omega_0 : c(x_i) = 1, c(y_j) = 0, 1 \le i \le m, 1 \le j \le n
\}. \\
\end{split}
\end{equation*}
where
$\{x_1, \dots ,x_m,y_1, \dots , y_n \} \subset \G$ and $m,n \ge 1$.
Each set in $\mathcal C$ has finite measure
since $\mu_0(S_x^0 \triangle S_y^0) < \infty$ for all $x,y \in \G$.
Assertion (a) follows from the fact that the family $\mathcal C$
is a countable base for the topology of $\Omega_0$
consisting of compact open sets.
Assertion (b) follows because every compact subset of $\Omega_0$
is contained in a finite union of sets in $\mathcal C$, and
hence is contained in a finite union of sets of finite measure.
\qed

\begin{proposition} \label{2}
The class of measure definite kernels on $\G$ is closed under pointwise
convergence.
\end{proposition}

{\sc Proof:}
By Proposition~\ref{1+} we may consider only measure definite kernels
of the form $\mu(S_x^0 \triangle S_y^0)$ for regular Borel measures
$\mu$ on $\Omega_0$.
Suppose that measures $(\mu_j)_{j=1}^{\infty}$ on $\Omega_0$ give rise
to measure definite kernels $f_j$ which converge pointwise to $f$.
We must show that $f$ is measure definite.
For fixed $x,y \in \G$, 
$f_j(x,y) \to f(x,y)$, so the  restricted measures
$$\mu_j |_{S_x^0 \triangle S_y^0}$$
are positive and of uniformly bounded finite mass.
Since the set of positive Borel measures of total mass
bounded by a given fixed constant form a compact metric space,
we can, by passing to a subsequence, assume that
$\mu_j |_{S_x^0 \triangle S_y^0}$
approaches some positive, finite mass Borel measure
on $S_x^0 \triangle S_y^0$.
Since $\G$ is countable, we may pass to a subsequence
where this happens for all $x,y \in \G$.

Write  $\Omega_0 = \bigsqcup _{k=1}^{\infty} T_k$, where $\mu_j | _{T_k}$
has a limit for each $T_k$ and where each $S_x^0 \triangle S_y^0$
is contained in a finite union of $T_k$'s.
To see that this is possible, 
suppose that the elements of $\G$ are labelled, $x_1$, $x_2$, etc.  Let
all sets of the two forms
$$
\{ c \in \Omega_0 :
   (c(x_1),\dots,c(x_{n-1}),c(x_n)) = (1,\dots,1,0) \}
$$
and
$$
\{ c \in \Omega_0 :
   (c(x_1),\dots,c(x_{n-1}),c(x_n)) = (0,\dots,0,1) \}
$$
make up the collection $(T_k)_k$.  This collection will be countable
or finite as $\G$ is countable or finite.  The disjoint union of the
$T_k$ gives all of $\Omega_0$.

 Put the limit measures 
on the $T_k$ together to form a measure $\mu$ on $\Omega_0$. For any
$x,y \in \G$, we have
\begin{equation*} 
\begin{split}
\mu(S_x^0 \triangle S_y^0)
&= \sum_{k=1}^{\infty}\lim_{j} \mu_j(T_k \cap (S_x^0 \triangle S_y^0)) 
= \lim_{j} \sum_{k=1}^{\infty}\mu_j(T_k \cap (S_x^0 \triangle S_y^0)) \\
&= \lim_{j}\mu_j (S_x^0 \triangle S_y^0) 
= \lim_{j} f_j(x,y) 
= f(x,y) . \\
\end{split}
\end{equation*}
This shows that $f$ is measure definite and so
proves the result.
\qed\

\medskip
Let $(u_j)_{j=1}^N \subset {\mathcal H}$ and $(u'_j)_{j=1}^N \subset {\mathcal H}'$
be two tuples of vectors in two real Hilbert spaces. We say that the two
tuples
are in the same {\it configuration\/} if 
$\|u_j-u_k\| = \|u'_j-u'_k\|$ for $1 \le j,k \le N$. 
If  ${\mathcal H}$,${\mathcal H}'$ have the same dimension, then the two tuples
are in the same configuration if and only if there is an affine isometry 
from ${\mathcal H}$ to ${\mathcal H}'$ sending $u_j$ to $u'_j$, $1 \le j \le N$. 

\begin{proposition} \label{3}
Let $f$ be a negative definite kernel on a set $\G$ and choose a mapping
$x \mapsto v_x$ of $\G$ into a real Hilbert space  ${\mathcal H}$
such that $f(x,y) = \|v_x - v_y\|^2$. Then

\noindent(i) $\sqrt f$ is measure definite;

\noindent(ii) there exists a regular Borel measure $\mu_0$ on $\Omega_0$
such that $f(x,y) =  \mu_0(S_x^0 \triangle S_y^0)$
for all $x,y \in \G$ and having following property:

if $\{x_1, \dots ,x_m,y_1, \dots , y_n \} \subset \G$, where $m,n \ge 1$,
then $\mu_0((S_{x_1}^0\cap \dots \cap S_{x_m}^0) \backslash
(S_{y_1}^0\cup \dots \cup S_{y_n}^0))$ is determined canonically by
the configuration of $(v_{x_1}, \dots ,v_{x_m},v_{y_1}, \dots , v_{y_n})$.
\end{proposition}

{\sc Proof:} Suppose the result has been proved for finite $\G$.  In
the infinite case, think of~$\G$ as the increasing union of a sequence
of finite sets, and for each of these finite sets construct a
measure~$\mu_0$ satisfying the conditions of the statement insofar as
they apply to elements of that set.  Take limits as in the proof of
Proposition~1.3 to obtain a measure~$\mu_0$ satisfying the conditions
for all of~$\G$.

Suppose, therefore, that $\G$~is finite.  Replace ${\mathcal H}$ by the
linear subspace which the the~$v_x$ generate.  This guarantees that
$\mathcal H$ is finite dimensional.  Let $\Omega$ be the homogeneous space
of all half-spaces of ${\mathcal H}$.  Then $\Omega$ may be realized as
$E_n/ E_{n-1}$, where $E_n$ is the (unimodular) group of rigid motions
of the $n$-dimensional real Hilbert space ${\mathcal H}$ and the
stabilizer of a given half space is isomorphic to $E_{n-1}$. Endow
$\Omega$ with its usual measure $\mu$, which is invariant under rigid
motions of ${\mathcal H}$.  We can normalize $\mu$ so that
$$\mu\{M \in \Omega : v \in M, w \notin M \} = \|v - w \| .$$
(This follows, for example, from the explicit expression for the
Haar measure on the affine group of ${\mathcal H}$ given in 
\cite[Chap.VII \S3 No.3 Ex.4]{bou}. It may also be seen geometrically.)
Let $S_x =  \{M \in \Omega : v_x \in M \}$.
Then ${\sqrt f}(x,y) = \|v_x - v_y\| = \mu(S_x \backslash S_y) 
= \mu(S_y \backslash S_x)$ is measure definite. This proves (i).

To prove (ii) it suffices to show that the measure
$$
\mu\{M \in \Omega : v_{x_j} \in M, v_{y_k} \notin M, 1 \le j \le m, 1 \le k
\le n\}
$$
depends only on the configuration of $(v_{x_1}, \dots ,v_{x_m},v_{y_1},
\dots , v_{y_n})$.
This is clear except for a possible dependence on 
the dimension of  ${\mathcal H}$.
Suppose then that 
$$
\{v_{x_1}, \dots ,v_{x_m},v_{y_1}, \dots , v_{y_n}\} \subset {\mathcal H}
\subset {\mathcal H}'
$$
Let $\Omega'$ be the space of all half-spaces of ${\mathcal H}'$ and let $\mu'$
be the natural 
measure on $\Omega'$. We must prove that 
\begin{multline*}
\mu\{M \in \Omega : v_{x_j} \in M, v_{y_k} \notin M, 1 \le j \le m, 1 \le k
\le n\} \\
= \mu'\{M' \in \Omega' : v_{x_j} \in M', v_{y_k} \notin M', 1 \le j \le m,
1 \le k \le n\} .
\end{multline*}
\noindent Define a measure $\nu$ on $\Om$ by
$\nu(E) = \mu'\{M' \in \Om' : M' \cap {\mathcal H} \in E \}$.
Since $\nu$ is invariant under rigid motions of ${\mathcal H}$ and satisfies
$$
\nu\{M \in \Om : v \in M, w \notin M \} = \mu'\{M' \in \Om' : v \in M', w
\notin M' \}
= \|v - w \|
$$
it follows that $\nu = \mu$. 
\qed
\medskip

\noindent {\bf Remarks}. (i) Let $\G$ be the set of vertices of a tree. The
usual
integer-valued distance function $d(x,y)$ between elements of $\G$ is a
measure
definite kernel. To see this, fix a vertex $o \in \G$ and for each $x \in
\G$
let $S_x = [o,x]$, the geodesic path between $o$ and $x$. Then $d(x,y)
=\mu(S_x \triangle S_y)$,
where $\mu$ is Lebesgue measure on the geometric realization of the tree.
(However it is easy to see that not all measure definite kernels arise
from the distance function on a tree.)
A. Valette \cite [Section 3]{v} has given an appropriate construction of
Lebesgue measure for real trees from which it follows similarly that the distance
function in a real tree is measure definite.

\begin{figure}[htbp]
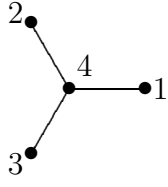
\label{tree}
{}\hfill
\beginpicture
\setcoordinatesystem units <0.5cm, .866cm>  
\setplotarea  x from -2.5 to 2.5,  y from -1.5 to 1.5
\put {$\bullet$} at 0 0
\put {$\bullet$} at -1 1
\put {$\bullet$} at -1 -1
\put {$\bullet$} at 2 0
\put {$4$} [b,l] at 0.2 0.2
\put {$3$} [r,t] at  -1.2 -1
\put {$1$} [l] at  2.2 0
\put {$2$} [b,r] at  -1.2  1
\putrule from 0 0 to 2 0
\setlinear \plot  -1 -1  0 0 -1 1 /
\endpicture
\hfill{}
\caption{A tree with four vertices}
\end{figure}

(ii) It is not true that $f^2$ is negative definite
whenever $f$ is measure definite. Let $\G = \{1,2,3,4\}$ be the set
 of vertices of the tree in Figure \ref{tree}.

\refstepcounter{picture}

\noindent Thus $d(1,4) = d(2,4) = d(3,4) = 1$ and $d(1,2) = d(2,3) = d(3,1)
= 2$
By the previous remark, $d(x,y)$ is measure definite.
There is clearly no set $\{v_1,v_2,v_3,v_4\}$  of vectors in a Hilbert
space
such that $d(x,y) = \|v_x - v_y\|$.
Therefore $d^2$ is not  negative definite. 

\medskip

\section{property (T)}

We  give an application of the preceding ideas to characterize 
countable groups with property (T)
in terms of  measure preserving actions on infinite measure spaces.
Another characterization has been given by K. Schmidt \cite[Proposition
2.10]{s},
A. Connes and B. Weiss \cite{cw}
in terms of  measure preserving actions on finite measure spaces,
using the notion of strong ergodicity.

We recall the Delorme--Guichardet characterization of groups with
property (T) \cite[Chapter~5, Theorem 20]{hav}, which may be expressed
as follows: a group has property (T) if and only if every left
$\G$-invariant negative definite kernel on $\G$ is bounded.

\begin{theorem} \label{4} A countable group $\G$ has property (T)
if and only if it satisfies the following condition:
for every measure-preserving action of $\G$  on a measure
space  $(\Omega,{\mathcal B},\mu)$ and every set  $S \in {\mathcal B}$  such that
$\mu (S \triangle gS) < \infty $  for all  $g \in \G$, we have  
 $$\sup_{g\in \G} \mu(S \triangle gS) < \infty .$$
\end{theorem}

{\sc Proof:}
Suppose that $\G$ has property (T), and that we are given a measure 
preserving action as described. Define a
left $\G$-invariant measure definite kernel on
$\G$ by $f(x,y) = \mu(xS \triangle yS)$.
By Proposition~\ref{1}, $f$ is negative definite and so
uniformly bounded, by the Delorme--Guichardet result.

Conversely, suppose that the condition in the theorem holds, and let
$f$ be a left $\G$-invariant negative definite kernel on $\G \times
\G$. By Proposition~\ref{3} $\sqrt f$ is measure definite with
corresponding regular Borel measure $\mu_0$ defined on the space
$\Omega_0 = \{0,1\}^{\G} \backslash \{(0,0,0, \dots ), (1,1,1, \dots
)\}$.  We must show that $\mu_0$ is invariant under the natural action
of $\G$ on $\Omega_0$ defined by $(gc)(x_0) = c(g^{-1}x_0)$.  Recall
that $\mathcal C$ is the family of all cylinder sets of the form
\begin{equation*} 
\begin{split}
A
&=  (S_{x_1}^0\cap \dots \cap S_{x_m}^0) \backslash
(S_{y_1}^0\cup \dots \cup S_{y_n}^0)\\
&= \{ c \in \Omega_0 : c(x_i) = 1, c(y_j) = 0, 1 \le i \le m, 1 \le j \le n
\} \\
\end{split}
\end{equation*}
where
$\{x_1, \dots ,x_m,y_1, \dots , y_n \} \subset \G$ and $m,n \ge 1$.
\noindent Observe that $A$ is compact and open in $\Omega_0$.

To pass from $A$ to $gA$, one must apply $g$ to each of the
$x_i$'s and $y_j$'s. Since
$\|v_{gx} - v_{gy}\|^2 = f(gx,gy) = f(x,y) = \|v_x - v_y\|^2$,
it follows from Proposition~\ref{3}(ii) that
the measures $\mu_0(A)$ and $\mu_0(g(A))$ coincide.
We have therefore shown that $\mu_0$ and $\mu_0^g$
coincide on all sets in ${\mathcal C}$.

Suppose that the elements of $G$ are labelled, $x_1$, $x_2$, etc.
Let $\mathcal C_1$ be the subset of $\mathcal C$
containing those sets of the form
$$
\{ c \in \Omega_0 : (c(x_1),c(x_2),\dots,c(x_n)) =
   (\epsilon_1,\epsilon_2,\dots,\epsilon_n) \}
$$
where $n$ is a positive integer, where each $\epsilon_j$ is either $0$
or $1$, and where each of the values $0$ and $1$ occurs at least once.

The sets in $\mathcal C_1$ satisfy the following properties: (1) $\mathcal
C_1$ is a neigborhood base for $\Omega$; (2) if two sets in $\mathcal C_1$
are not disjoint, then one contains the other; (3) $\mathcal C_1$ contains
no infinite ascending chain.


Let $U$ be an open subset of $\Omega$.  It follows from the three
properties of $\mathcal C_1$ that $U$ is the disjoint union of those sets in
${\mathcal C_1}$ maximal with respect to being contained in $U$.
%
%
Thus $\mu_0(U) = \mu_0(g(U))$.  Since the measures $\mu_0$ and
$\mu_0^g$ are regular and they coincide on open sets, they must agree
on all Borel sets.  Therefore the action of $\G$ on $\Omega_0$ is
measure-preserving.


Now write $S = S_e^0$. Since $\sqrt f$ is measure definite,
$\mu_0(S \triangle gS) < \infty $  for all  $g \in \G$.
The condition in the theorem implies that
$\sup_{g\in \G} \mu_0(S \triangle gS) < \infty$.
The negative definite kernel $f$ is therefore bounded
and so $\G$ has property (T).
\qed 

\medskip
\noindent {\bf Questions}. Several questions remain.

(i) Can measure definite kernels be characterized internally, without
reference to measure spaces?

(ii) Is there a result analogous to Theorem \ref{4} for locally compact
groups?

(iii) Can the action of $\G$ in Theorem \ref{4} be assumed to be ergodic?


\begin{thebibliography}{RRR}

\bibitem [AW]{akw} C.A. Akemann and M.E. Walter, Unbounded negative
definite functions,
{\it Canad. J. Math. }{\bf33 }(1981), 862-871.

\bibitem[B]{bou} N.\ Bourbaki, {\em Int\'egration, Chapitres VII et VIII}, 
Hermann, Paris, 1963.

\bibitem[CW]{cw} A. Connes and B. Weiss, Property (T) and asymptotically
invariant
sequences, {\it Israel J. Math. }{\bf37} (1980), 209-210.

\bibitem [D]{del} P. Delorme, 1-cohomologie des repr\'esentations unitaires
des groupes
de Lie semi-simples et r\'esolubles. Produits tensoriels et
repr\'esentations,
{\it Bull. Soc. Math. France }{\bf105} (1977), 281-336.

\bibitem [DR]{dr} A. Deutsch and G. Robertson, Functions of conditionally
negative type
 on Kazhdan groups, {\it Proc. Amer. Math. Soc. }{\bf123} (1995), 919-926.

\bibitem [G]{gui} A. Guichardet, \'Etude de la 1-cohomologie et la
topologie du dual
pour les groupes de Lie \`a radical ab\'elien, {\it Math. Ann. }{\bf 228}
(1977), 215-232.

\bibitem [HV]{hav} P. de la Harpe and A. Valette, La propri\'et\'e (T) de
Kazhdan pour 
les groupes localement compacts, {\it Ast\'erisque }{\bf 175}, Soc. Math.
France, 1989.

\bibitem [K]{k} D. Kazhdan, Connection of the dual space of a group with
the structure
 of its closed subgroups, {\it Funct. Anal. and its Appl.} {\bf 1} (1967),
63-65.

\bibitem [R]{r} W. Rudin, {\em Real and Complex Analysis}, 2nd. edition,
McGraw-Hill,
New-York, 1974.

\bibitem [S]{s} K. Schmidt, Asymptotically invariant sequences and an
action of
 $SL(2,\bold{Z})$ on the 2-sphere, {\it Israel J. Math.} {\bf 37} (1980),
193-208.

\bibitem [V]{v} A. Valette, Les repr\'esentations uniform\'ement born\'ees
associ\'ees
a un arbre r\'eel, Bull. Soc. Math. Belgique, {\bf152} (1990), 747-760.

\end{thebibliography}
\end{document}